\documentclass[12pt,psamsfonts,leqno,oneside,letterpaper]{amsart}
\usepackage[dvips,text={6.5truein,9truein},left=1truein,top=1truein]{geometry}
\usepackage{amssymb,amsmath,amscd,enumerate,
}
\usepackage[pdftex]{graphicx}
\usepackage{url}

\usepackage[colorlinks,linkcolor=blue,citecolor=blue,pdfstartview=FitH]{hyperref}
\input xy
\xyoption{all}
\SelectTips{cm}{12}

\usepackage{color}
\usepackage{mathrsfs}

\newcommand\missingref[1]{\empty}
\newcommand\tinymissingref[1]{\empty}
\newcommand\abstractcomment[1]{\empty}

\theoremstyle{definition}
\newtheorem{para}{}[section]

\newtheorem{remark}[para]{Remark}
\newtheorem{remarks}[para]{Remarks}
\newtheorem{notation}[para]{Notation}

\newtheorem{convention}[para]{Convention}
\newtheorem{definition}[para]{Definition}
\newtheorem{definitions}[para]{Definitions}
\newtheorem{definitionnotation}[para]{Definition and Notation}
\newtheorem{remarksnotation}[para]{Remarks and Notation}

\newtheorem{remarknotation}[para]{Remark and Notation}
\newtheorem{notationremark}[para]{Notation and Remark}
\newtheorem{notationremarks}[para]{Notation and Remarks}
\newtheorem{definitionsnotation}[para]{Definitions and Notation}

\newtheorem{reviewdefinition}[para]{Review and Definition}

\newcommand\Alternatives{\begin{enumerate}[(i)]}
\newcommand\EndAlternatives{\end{enumerate}}
\newcommand\Conditions{\begin{enumerate}[(1)]}
\newcommand\EndConditions{\end{enumerate}}

\theoremstyle{plain}
\newtheorem{theorem}[para]{Theorem}
\newtheorem{lemma}[para]{Lemma}
\newtheorem{proposition}[para]{Proposition}
\newtheorem{maintheorem}[para]{Main Theorem}
\newtheorem{corollary}[para]{Corollary}
\newtheorem{conjecture}[para]{Conjecture}

\newtheorem{claim}[equation]{}
\numberwithin{equation}{para}
\numberwithin{figure}{section}

\newcommand\Number{\begin{para}}
\newcommand\EndNumber{\end{para}}
\newcommand\Definition{\begin{definition}}
\newcommand\EndDefinition{\end{definition}}
\newcommand\Definitions{\begin{definitions}}
\newcommand\DefinitionsNotation{\begin{definitionsnotation}}
\newcommand\RemarksNotation{\begin{remarksnotation}}
\newcommand\EndRemarksNotation{\end{remarksnotation}}
\newcommand\RemarkNotation{\begin{remarknotation}}
\newcommand\EndRemarkNotation{\end{remarknotation}}
\newcommand\NotationRemark{\begin{notationremark}}
\newcommand\NotationRemarks{\begin{notationremarks}}
\newcommand\EndNotationRemark{\end{notationremark}}
\newcommand\EndNotationRemarks{\end{notationremarks}}
\newcommand\EndDefinitionsNotation{\end{definitionsnotation}}
\newcommand\DefinitionNotation{\begin{definitionnotation}}
\newcommand\EndDefinitionNotation{\end{definitionnotation}}
\newcommand\ReviewDefinition{\begin{reviewdefinition}}
\newcommand\EndReviewDefinition{\end{reviewdefinition}}
\newcommand\EndDefinitions{\end{definitions}}
\newcommand\Theorem{\begin{theorem}}
\newcommand\EndTheorem{\end{theorem}}
\newcommand\Conjecture{\begin{conjecture}}
\newcommand\EndConjecture{\end{conjecture}}
\newcommand\Remark{\begin{remark}}
\newcommand\EndRemark{\end{remark}}
\newcommand\Remarks{\begin{remarks}}
\newcommand\EndRemarks{\end{remarks}}
\newcommand\Convention{\begin{convention}}
\newcommand\EndConvention{\end{convention}}
\newcommand\Notation{\begin{notation}}
\newcommand\EndNotation{\end{notation}}
\newcommand\Lemma{\begin{lemma}}
\newcommand\EndLemma{\end{lemma}}
\newcommand\Proposition{\begin{proposition}}
\newcommand\EndProposition{\end{proposition}}
\newcommand\Corollary{\begin{corollary}}
\newcommand\EndCorollary{\end{corollary}}
\newcommand\Claim{\begin{claim}}
\newcommand\EndClaim{\end{claim}}
\newcommand\Proof{\begin{proof}}
\newcommand\EndProof{\end{proof}}
\newcommand\Equation{\begin{equation}}
\newcommand\EndEquation{\end{equation}}

\newcommand\Bullets{\begin{itemize}}
\newcommand\EndBullets{\end{itemize}}

\renewcommand\epsilon{\varepsilon}

\newcommand\tP{\widetilde P}

\newcommand\tS{{\widetilde S}}

\newcommand\calc{{\mathcal C}}

\newcommand\calf{{\mathcal F}}

\newcommand\calg{{\mathcal G}}

\newcommand\calu{{\mathcal U}}

\newcommand\caly{{\mathcal Y}}

\newcommand\sigmabar{{\overline{\sigma}}}
\newcommand\calw{{\mathcal W}}
\newcommand\calv{{\mathcal V}}

\newcommand\Isom{{\rm Isom}}

\newcommand\ZZ{{\bf Z}}

\newcommand\HH{{\bf H}}

\newcommand\iccg{ICC-group}
\newcommand\link{\mathop{\rm link}}

\newcommand\calz{{\mathcal Z}}

\newcommand\rank{\mathop{{\rm rank}}}
\newcommand\shmel[2]{{#1}_{\tiny{{\rm IR}\le {#2}}}}

\DeclareFontFamily{U}{rcjhbltx}{}
\DeclareFontShape{U}{rcjhbltx}{m}{n}{<->rcjhbltx}{}
\DeclareSymbolFont{hebrewletters}{U}{rcjhbltx}{m}{n}

\let\aleph\relax\let\beth\relax
\let\gimel\relax\let\daleth\relax

\DeclareMathSymbol{\aleph}{\mathord}{hebrewletters}{39}
\DeclareMathSymbol{\beth}{\mathord}{hebrewletters}{98}
\DeclareMathSymbol{\gimel}{\mathord}{hebrewletters}{103}
\DeclareMathSymbol{\daleth}{\mathord}{hebrewletters}{100}

\DeclareMathSymbol{\lamed}{\mathord}{hebrewletters}{108}
\DeclareMathSymbol{\mem}{\mathord}{hebrewletters}{109}
\DeclareMathSymbol{\ayin}{\mathord}{hebrewletters}{96}
\DeclareMathSymbol{\tsadi}{\mathord}{hebrewletters}{118}
\DeclareMathSymbol{\qof}{\mathord}{hebrewletters}{114}
\DeclareMathSymbol{\shin}{\mathord}{hebrewletters}{152}

\title{The geometry of $\MakeLowercase{k}$-free hyperbolic $3$-manifolds}
\author{Rosemary K. Guzman}
\address{Department of Mathematics
\\
University of Chicago\\
5734 S. University Ave.\\
Chicago, IL 60637-1514}
\email{guzman@math.uchicago.edu}

\author{Peter B. Shalen}
\address{Department of Mathematics, Statistics, and Computer Science
(M/C 249)\\
University of Illinois at Chicago\\
851 S. Morgan St.\\
Chicago, IL 60607-7045}
\email{petershalen@gmail.com}

\begin{document}

\maketitle
\begin{abstract}
We investigate the geometry of closed, orientable, hyperbolic $3$-manifolds whose fundamental groups are $k$-free for a given integer $k\ge 3$.  We show that any such manifold $M$ contains a point $P$  with the following property:
If $S$ is the set of elements of $\pi_1(M,P)$ represented by loops of
length $<\log(2k-1)$, then for every subset $T \subset S$, we have
$\rank\langle  T\rangle \le k-3$. This generalizes to all $k\ge3$
results proved in \cite{ACCS} and \cite{cs}, which have been used to relate the volume of a hyperbolic manifold to its topological properties, and it strictly improves on the result obtained in \cite{Guzman} for $k=5$. The proof avoids the use of results about ranks of joins and intersections in free groups that were used in \cite{cs} and \cite{Guzman}.
    \end{abstract}

\section{Introduction}
{A group $\Gamma$ is said to be {\it $k$-free}, where $k$ is a positive integer, if every subgroup of $\Gamma$ having rank at most $k$ is free. We will say that a hyperbolic manifold $M$ is {\it $k$-free} if $\pi_1(M)$ is $k$-free.} 

 In \cite{Paradoxical}, \cite{ACCS}, \cite{ACS-surgery}, \cite{ACS}, \cite{CS_vol},  \cite{cs} and \cite{Guzman}, interesting and novel connections were drawn between { the geometry of a closed, orientable hyperbolic $3$-manifold $M$ and its topological properties. These results may be thought of as making Mostow rigidity more explicit. The topological property that arose in some of the main results of
\cite{Paradoxical}, \cite{ACCS},  \cite{cs} and \cite{Guzman} was the property of being $k$-free for a suitable positive integer $k$. The results of \cite{ACS}, \cite{CS_vol} are deduced from the results of \cite{ACS-surgery}, \cite{ACS} and \cite{CS_vol} via topological results that relate $k$-freeness, for a given  value of $k$, to more familiar topological information about $M$, such as the dimension of the $\ZZ_p$-vector space $H_1(M;\ZZ_p)$ for a suitable prime $p$, or the presence of a low-genus incompressible surface. Some of the results in 
\cite{ACCS}, \cite{ACS-surgery}, \cite{ACS}, \cite{CS_vol} and  \cite{cs} give strikingly strong upper bounds for the dimension of $H_1(M;\ZZ_p)$  in terms of upper bounds for the volume of $M$; these results are examples of ``quantitative Mostow rigidity.'' 

In order to 
relate quantitative geometric properties of $M$ to $k$-freeness one needs a connection between quantitative hyperbolic geometry and freeness of groups. The fundamental tool used in doing this is the ``$\log(2k-1)$ Theorem'' that was established in \cite{Paradoxical}, \cite{ACCS} and \cite{ACS-surgery}. The theorem asserts that if $x_1,\ldots,x_k$ are elements of $\Isom_+(\HH^3)$ that generate a free, discrete group of rank $k$, then for every point $\tP\in\HH^3$ we have $\sum_{i=1}^k1/(1+\exp(d(\tP,x_i\cdot \tP)))\le1/2$; in particular $\max_{1\le i\le k}d(\tP,x_i\cdot \tP)\ge\log(2k-1)$.  (The version of the theorem that we have quoted is \cite[Theorem 4.1]{ACS-surgery}. In the body of the paper we quote only a weaker version, for the case when $\langle x_1,\ldots,x_k\rangle$ is purely loxodromic; this follows immediately from the results of \cite{ACCS} and \cite{Agol_Tameness} or \cite{Calegari_Gabai}, and is used via its consequence \cite[Lemma 2.5]{Guzman}.)

In this paper, by refining the machinery used in some of the papers mentioned above, we have obtained stronger results about the connections between geometry and $k$-freeness than the ones given in those papers. Remarkably, we have done this while avoiding the difficult results about free groups that were used in the earlier arguments. Here is our main  result:

\begin{maintheorem}\label{Main Theorem}  Suppose that $M$ is a {$k$-free,} closed, orientable, hyperbolic $3$-manifold,  and that $k\ge 3$ is an integer.
Then there exists a point $P$ of $M$ with the following property.
If $S$ is the set of elements of $\pi_1(M,P)$ represented by loops of length $<\log(2k-1)$, then for every subset $T \subset S$ we have $\rank\langle  T \rangle\le k-3$.
\end{maintheorem}

For $k=3$, this main theorem is equivalent to the result, first proved in \cite{ACCS}, that if (a closed, orientable, hyperbolic $3$-manifold) $M$ is $3$-free, then $M$ contains (an isometric copy of) a hyperbolic $3$-ball of radius $(\log5)/2$; this can be expressed by saying that $M$ contains a ``$(\log5)$-thick''
 point. This fact was used in \cite{ACCS} to show that if $M$ is $3$-free then its volume is strictly bounded below by $3.08$. The latter fact was used in \cite{CS_vol} (building on topological work in \cite{ACS} and on geometric results from \cite{AST}) to show that if 
$M$ has volume at most $3.08$, then the $\ZZ_2$-vector space $H_1(M;\ZZ_2)$ has dimension at most $5$. 

For $k=4$, our main theorem is equivalent to Theorem 1.4 of \cite{cs}, which in the language of that paper says that if $M$ is $4$-free then it contains a ``$(\log7)$-semithick'' point. This fact is used in \cite{cs} to show that if $M$ is $4$-free then its volume is strictly bounded below by $3.44$. The latter fact was used in \cite{cs},  in combination  with the results of \cite{AST} and \cite{CS_vol}, to show that if 
$M$ has volume at most $3.44$, then $H_1(M;\ZZ_2)$ has dimension at most $7$.

For an arbitrary closed, orientable, hyperbolic $3$-manifold $M$, the best known general upper bound for the dimension of $H_1(M;\ZZ_2)$  in terms of the volume of $M$, given by
\cite[Proposition 2.2]{AgolLeiningerMargalit}, is $ 334.08$ times the volume. Note that the estimates established in 
 \cite{CS_vol} and \cite{cs} improve on this by a couple of orders of magnitude when the volume is suitably small.

In the case  $k=5$, our main theorem is strictly stronger than the main theorem of \cite{Guzman}; whereas in this case our main theorem asserts that every subset of the set $S$ generates a subgroup of rank at most  $2$, the result established in \cite{Guzman} asserts only that $S$ is contained in some subgroup of rank at most $2$.

{
Because our main theorem generalizes the results of \cite{ACCS} and \cite{cs}, which were used to give good lower bounds for volumes of $3$-free and $4$-free hyperbolic manifolds, it is very natural to try to use them to give good lower bounds for volumes of $k$-free hyperbolic manifolds, and we  hope to do this in a future paper.
The method of proof of  \cite[Proposition 2.2]{AgolLeiningerMargalit} can be expected to give a lower bound of $\alpha\cdot k$ for the volume of a $k$-free, closed, orientable hyperbolic $3$-manifold, where $\alpha$ is of the order of magnitude of $100$; the results of \cite{ACCS} and \cite{cs} suggest the possibility of replacing $\alpha$ here by a constant of the order of magnitude of $1$, and this will be our goal.

Our plan for using Main Theorem \ref{Main Theorem} to obtain a good
lower bound for the volume of  $M$ involves
extending the arguments used in \cite{cs} for the case $k=4$. If we
write $M=\HH^3/\Gamma$, where $\Gamma\le\Isom_+(\HH^3)$ is discrete,
cocompact and $k$-free, the conclusion of Main Theorem \ref{Main
Theorem} may be reformulated as saying that there is a point
$\tP\in\HH^3$ such that, if $\tS$ denotes the set of elements
$\gamma\in\Gamma$ such that $d(\tP,\gamma\cdot\tP)<\log(2k-1)$, then
the subgroup of $\Gamma$ generated by any subset of $\tS $ has rank at
most $k-3$. In the case $k=4$ this means that $\tS $ is contained in a
cyclic subgroup $C$ of $\Gamma$, and in \cite{cs} it is shown that
unless $M$ contains a very short geodesic---a case that must be
treated separately---the possibilities for the set $\tS \subset C$ are
strongly limited. This control over the set $\tS $ is used to give
geometric restrictions on the configuration of points $\tS \cdot
\tP=\{\gamma\cdot \tP:\gamma\in \tS
\}$, which are in turn used to give a lower bound for the volume of the
image in $M$ of the ball of radius $\log(2k-1)/2$ centered at $\tP$,
and hence for the volume of $M$.

For $k>4$, our preliminary investigations indicate that the possibilities
for the set $\tS$ are severely limited by the condition that every
subset of $\tS$ generates a subgroup of rank at most $m:=k-3$. Indeed,
$F:=\langle \tS\rangle$ has rank at most $m$, and is free since $\Gamma$
is $k$-free. It appears that in a free group $F$ of rank $m$ (say),
the subsets $\tS$ with the property that every subgroup of $\tS$ generates
a subgroup of rank at most $m$ are so special that one can expect to
classify them up to automorphisms of $F$. Our preliminary calculations
for particular examples of the set $\tS$ also suggest that this will
give a lower bound for the volume of the desired order of magnitude.

Eventually we  hope to be able to use the projected results about volumes of $k$-free manifolds to improve the constant in
 \cite[Proposition 2.2]{AgolLeiningerMargalit} by a couple of orders of magnitude, but this will require extending in a non-trivial way the arguments given in \cite{CS_vol} and \cite{cs} for  bounding the dimension of $H_1(M;\ZZ_2)$ in the case where $M$ contains a low-genus incompressible surface.
}


{The proof of our main theorem will be given in Section 5, building on results and observations in Sections 2--4. For the proof} we write $M=\HH^3/\Gamma$, where $\Gamma$ is a $k$-free (and hence torsion-free), cocompact discrete subgroup of $\Isom_+(\HH^3)$. For each maximal cyclic subgroup $C$ of $\Gamma$ generated by an element of translation length less than $\lambda:=\log(2k-1)$, we have a ``hyperbolic cylinder'' $Z_\lambda(C)$ consisting of all points $\tP\in\HH^3$ such that $d(\tP,\gamma\cdot\tP)<\lambda$ for some nontrivial element $\gamma$ of $C$. We regard the $Z_\lambda(C)$ as forming a family of hyperbolic cylinders, indexed by maximal cyclic subgroups each of which is generated by an element of translation length less than $\lambda$; if the conclusion of the theorem is false, it is easily deduced that this family covers $\HH^3$, and the context for the proof is the study of the (geometrically realized) nerve $K$ of this covering. (The definitions of the nerve of a covering, and of other terms used in this sketch, are reviewed in Section \ref{prelim section}.) Since the covering is indexed by certain maximal cyclic subgroups of $\Gamma$, the complex $K$ comes equipped with a ``labeling'' which assigns to each vertex $v$ a maximal cyclic subgroup $C_v$ of $\Gamma$. If for each (open) simplex $\tau$ of $K$ we denote by $\Theta(\tau)$ the  subgroup of $\Gamma$ generated by the cyclic subgroups $C_v$ where $v$ ranges over the vertices of $\tau$, we define the `` internal rank'' of a simplex $\sigma$ to be the maximum of the ranks of $\Theta(\tau)$, where $\tau$ ranges over all faces of $\sigma$, including $\sigma$ itself (see Definitions \ref{nr} and \ref{more Theta}). The $\log(2k-1)$ Theorem (via its consequence \cite[Lemma 2.5]{Guzman}), with the hypothesis that $\Gamma$ is $k$-free, is used to show that each simplex of $K$ has internal rank at most $k-1$.  (It is at this step that the choice of the value $\lambda=\log(2k-1)$ is crucial.)

For every positive integer $m$, the simplices of $K$ having internal rank at most $m$ form a subcomplex of $K$, which we denote 
by {$\shmel Km$}. If the conclusion of the main theorem is assumed to be false, we show, using an argument based on the Borsuk Nerve Theorem, the elementary topology of simplicial complexes, and convexity considerations in $\HH^3$, that the subset {$|K|-|\shmel K{k-3}|$} of the underlying space $|K|$ of $K$ is contractible. This is the content of Lemma \ref{ContractibilityLemma}, which is the main result of Section \ref{contractibility section}.

Since every simplex of $K$ has internal rank at most $k-1$, we may express {$|K|-|\shmel K{k-3}|$} as a set-theoretical disjoint union of $X_{k-2}$ and $X_{k-1}$, where $X_r$ denotes the union of all simplices of internal rank $r$ for $r=k-2,k-1$. Each set $X_r$ is saturated in the sense that it is a union of (open) simplices. Generalizing the definition of $\Theta(\tau)$ given above, we define $\Theta(X)$, for every saturated subset $X$ of $K$, to be the subgroup of $\Gamma$ generated by the cyclic groups $C_v$, where $v$ ranges over all vertices of simplices in $X$. A combinatorial argument about subgroups of free groups, making strong use of the hypothesis that $\Gamma$ is $k$-free, shows that  for $r=k-2,k-1$ the group $\Theta(X_r)$  is locally free, in the sense that its finitely generated subgroups are all free; furthermore, the local rank of $X_r$, defined to be the supremum of the  ranks of its finitely generated subgroups, is at most $r$. This is the content of Lemma \ref{RankLemma}, which is the main result of Section \ref{rank section}.

In Section \ref{main proof section} we describe a bipartite graph $T$ whose vertices correspond to connected components of the sets $X_r$ for $r=k-2,k-1$. The contractibility of {$|K|-|\shmel K{k-3}|$} is used to show that $T$ is a tree. Naturality considerations (encoded in the notion of a ``label-preserving action'' which is defined in Section \ref{prelim section}) give rise to an action, without inversions, of $\Gamma$ on $T$. By combining the local freeness of the groups $\Theta(X_{k-2})$ and $\Theta(X_{k-1})$, and the bounds on their local ranks, with further combinatorial arguments about free groups and the $k$-freeness of $\Gamma$, we deduce that the stabilizers of vertices under this action are locally free. But a $3$-manifold topology argument shows that the fundamental group of a closed aspherical manifold cannot admit an action without inversions on a tree under which all vertex stabilizers are locally free. This contradiction completes the proof of the main theorem.

{This method of proof has much in common with the method used in \cite{cs} and \cite{Guzman}. One crucial new  feature is the use of the notion of internal rank of a simplex $\sigma$. In \cite{cs} the arguments involved the rank of $\Theta(\sigma)$, while in \cite{Guzman} they involved the ``minimum enveloping rank'' of $\Theta(\sigma)$, which is defined to be the minimum rank of any subgroup of $\Gamma$ containing $\Theta(\sigma)$. One advantage of using internal rank in place of minimum enveloping rank is that the conclusion of the  main theorem of the present paper is stronger than the main result of \cite{Guzman} even for $k=5$, as noted above. 

The other advantage of using internal rank in place of minimum enveloping rank is that it allows us, in the combinatorial argument required for the proof of Lemma \ref{RankLemma} of this paper, to avoid the difficult group-theoretical issues that arise in the proofs of the analogous results
  \cite[Proposition 4.4]{cs} and \cite[Proposition 4.3]{Guzman}. Indeed, the arguments needed to prove the latter results involve
the following group-theoretical conjecture, the general case of which is formulated in \cite{Guzman}: if $m\ge1$  is an integer, and if $E$ and $E'$ are rank-$m$ subgroups of a free group $F$ such that $\rank(E\cap E')\ge m$, then $\rank\langle E,E'\rangle\le m$. This was proved in \cite{Kent} and \cite{LMc} for $m=2$, and in \cite{Guzman} for $m=3$. The arguments of \cite{Guzman} show that if $k\ge4$ is an integer such that if this conjecture is true for $m=k-2$, then the main result of \cite{Guzman} remains true when ``$5$-free'' is replaced in the hypothesis by ``$k$-free.'' However, it was shown in \cite{SourGrapes} that this group-theoretical conjecture is false for $m\ge5$. Thus when $k\ge7$, the use of internal rank seems indispensable even for getting results similar to the main result of \cite{Guzman} with ``$k$-free'' in place of ``$5$-free''; and even when $k=5$, the use of internal rank gives a strictly stronger result.

We are indebted to Marc Culler and Benson Farb for valuable discussions of the material presented in this paper, and for their encouragement.}

\section{Preliminaries}\label{prelim section}


\Definition A group will be said to have the {\it infinite cyclic centralizer property}, or more briefly to be an {\it \iccg}, if the centralizer of every non-trivial element of $\Gamma$ is infinite cyclic. Note that if $\Gamma$ is an \iccg\ then $\Gamma$ is torsion-free, and
 every non-trivial element $\gamma$ of $\Gamma$ belongs to a unique maximal cyclic subgroup of
$\Gamma$, namely the centralizer of $\gamma$ in $\Gamma$.
\EndDefinition

\Definition \label{cylinders} If $\Gamma$ is an \iccg, we will denote by $\calc(\Gamma)$ the set of maximal cyclic subgroups of $\Gamma$. \EndDefinition

\Definition \label{nr} If $\Gamma$ is an \iccg\ and $S$ is a finite subset of ${\mathcal C}(\Gamma)$, we define the {\it internal rank} of $S$, denoted ${\rm IR}(S)$, as the maximum of the set $\{ {\rm rank} \langle T \rangle : T \subseteq S \}$.\EndDefinition

\Number\label{complex stuff}
Except when we specify otherwise, the term {\it simplicial complex} will be understood in the geometric sense; that is, a simplicial complex is a set $L$  of pairwise-disjoint finite-dimensional open simplices in a (possibly infinite-dimensional) real vector space, with the property that any face of any simplex in $L$ is itself in $L$. The geometric realization of an abstract simplicial complex is a simplicial complex in this sense, and every simplicial complex is simplicially isomorphic to the geometric realization of an abstract simplicial complex.  The simplicial complexes referred to in this paper are not assumed to be locally finite. If $L$ is a simplicial complex, the union of its simplices will be denoted by $|L|$. The set $|L|$ will always be understood to be endowed with the weakest topology which induces the standard topology on the closure of each simplex. We will often denote the vertex set of a simplicial complex $L$ by $L^{(0)}$.

If $\sigma$ is a simplex of a simplicial complex $L$, we will denote by $\sigmabar$ the subcomplex of $L$ consisting of $\sigma$ and its faces (so that the topological closure of $\sigma$ in $|K|$ is $|\sigmabar|$). Thus $\sigmabar^{(0)}$ is the set of vertices of $\sigma$.

\EndNumber

\Definition Let  $L$ be a simplicial complex. We say that a subset $W$  of   $|L|$  is {\it saturated}  if $W$ is a union of (open)
simplices of $L$.
\EndDefinition

\Number\label{first nerve}
Recall that an indexed family $\calf=(U_i)_{i\in I}$ of nonempty (open) subsets of a topological space $X$ is said to {\it cover} $X$ (or to be a(n open) {\it covering} of $X$) if $X=\bigcup_{i\in I}U_i$. (Here the index set $I$ can be any set whatsoever.) We define the {{\it abstract nerve}} of a covering $\calf=(U_i)_{i\in I}$ of $X$ to be the abstract simplicial complex that is well defined up to canonical simplicial isomorphism as follows. The vertex set $V$ of the complex is a bijective copy of the index set $I$, equipped with a specific bijection $i\mapsto v_i$ from $I$ to $V$. A simplex $\sigma$  is a set  $\{v_{i_0},\ldots,v_{i_d}\}$, with $d\ge0$ and $i_0,\cdots,i_d\in I$, such that
$U_{i_0}\cap\cdots\cap U_{i_d}\ne\emptyset$. 
{The  {\it nerve} of a covering is defined to be the geometric realization of its abstract nerve.}
The nerve of a covering $\calf$ will be denoted $K_\calf$.

(In the most classical definition, an open covering of $X$ is a collection $\calu$ of open subsets of $X$. The vertex set of the nerve is the set $\calu$, and a simplex is a finite nonempty set of elements of $\calu$ whose intersection is nonempty. If an indexed covering $\calf=(U_i)_{i\in I}$ is ``faithfully indexed'' in the sense that for any two distinct indices $i,j\in I$ we have $U_i\ne U_j$, then our definition of the nerve is interchangeable with this more classical definition. However, in the proof of Lemma \ref{ContractibilityLemma}, we will need to consider indexed coverings which may not be faithfully indexed in this sense.
It is for this reason that the nerve construction for an indexed covering  is needed here.)


We will make essential use of the Borsuk Nerve Theorem \cite[Corollary 4G.3]{hatcher}: Let $\calf=(U_i)_{i\in I}$ be  an open
covering of a paracompact space $X$. Suppose that for every nonempty finite set $\Phi\subset I$ for which $W:=\bigcap_{i\in\Phi}U_i\ne\emptyset$, the set $W$ is contractible. Then $K_\calf$ is homotopy-equivalent to $X$.

\EndNumber

\Definitions {We review the following definitions from \cite{cs} and \cite{Guzman}.} 
Let $\Gamma$ be an \iccg.
We define {a {\it $\Gamma$-labeling} of a simplicial complex $L$ to be an indexed family $(C_v)_{v\in L^{(0)}}$, where  $C_v$ is a maximal cyclic subgroup of $\Gamma$ for every vertex $v$ of $L$. We define a {\it 
$\Gamma$-labeled complex} to be a pair $(L, (C_v)_{v\in L^{(0)}})$, where $L$ is a  simplicial complex and $(C_v)_{v\in L^{(0)}}$ is a labeling of $L$.}
We will often say that $L$ is a $\Gamma$-labeled complex when the family $(C_v)_{{v\in L^{(0)}}}$ is understood. By a   \emph{labeling-compatible action} of $\Gamma$ on {a $\Gamma$-labeled complex $(L, (C_v)_{{v\in L^{(0)}}})$,} we mean a simplicial action of $\Gamma$ on $L$ such that for every $\gamma\in\Gamma$ and every vertex $v$ of $L$  we have $C_{\gamma\cdot v}
= \gamma C_v \gamma ^{-1}$. 
\EndDefinitions

\Definitions\label{more Theta} Let $\Gamma$ be an \iccg, and $(L, (C_v)_{{v\in L^{(0)}}})$ a $\Gamma$-labeled complex. For any simplex $\sigma$ of $L$ we define $I_\sigma$ to be the subset  {$\{ C_v : v \in \sigmabar^{(0)}\} $} of $\Gamma$, and we define $\Theta(\sigma)$ to be the subgroup $
\langle I_\sigma\rangle=
 \langle C_v : v \in {\sigmabar^{(0)}} \rangle $ of $\Gamma$. 
We define the \emph{internal rank} of $\sigma$ by ${\rm IR} (\sigma) = {\rm IR} (I_ \sigma )$. In view of the definitions, it follows that ${\rm IR} (\sigma) = \max\{\rank\Theta(\tau) : \tau \le \sigma \}$.

These definitions may be extended to saturated sets as follows.
If $W$ is any saturated subset of $|L|$, we define a set $I_W\subset\Gamma$ by $I_W=\bigcup_\sigma I_\sigma$, where $\sigma$ ranges over all simplices of $L$ contained in $W$. We also set
$\Theta(W)=\langle I_W\rangle$, so that $\Theta(W)$ is the subgroup of $\Gamma$ generated by all subgroups of the form $\Theta(\sigma)$, where $\sigma$ ranges over all  simplices contained in $W$.

{Since ${\rm IR} (\sigma) = \max_{\tau \le \sigma} \rank\Theta(\tau)$ for every $\sigma\in L$, we have ${\rm IR}(\sigma')\le{\rm IR}(\sigma)$ whenever $\sigma$ is a simplex of $L$ and $\sigma'$ is a face of $\sigma$. Hence for every positive integer $m$, the set of all simplices of $K$ having internal rank at most $m$ is a subcomplex of $L$. It will be denoted by 
{$\shmel{L}{m}$.}
}
\EndDefinitions

{
\Number\label{why compatible}
Suppose that $\Gamma$ is an \iccg\ and that $(L, (C_v)_{{v\in L^{(0)}}})$ is a $\Gamma$-labeled complex equipped with a labeling-compatible action of $\Gamma$. If $\sigma$ is  any simplex of $L$, then since the action of $\Gamma$ is simplicial, we have 
$\overline{\gamma\cdot\sigma}^{(0)}=\gamma\cdot(\sigmabar^{(0)})$. The definition of a labeling-compatible action gives 
$C_{\gamma\cdot v}
= C_v^ \gamma$ for each $v\in \sigmabar^{(0)}$, where $x\mapsto x^\gamma$ denotes the inner automorphism $x\to\gamma x\gamma^{-1}$ of $\Gamma$. We therefore have $I_{\gamma\cdot\sigma}=\{C_{\gamma\cdot v}:v\in \sigmabar^{(0)}\}
=\{  C_v ^\gamma :v\in \sigmabar^{(0)}\}= I_\sigma^\gamma$. It follows that, more generally, if $W$ is a saturated subset of $L$, we have
$I_{\gamma\cdot W}=\bigcup_{\sigma\subset W}
I_{\gamma\cdot\sigma}=\bigcup_{\sigma\subset W} I_\sigma^\gamma=
 I_W^\gamma$.
 Hence $\Theta(\gamma\cdot W)=\langle I_{\gamma\cdot W}\rangle=\langle  I_W^\gamma\rangle= \Theta(W)^\gamma$. We have shown that: 
\Equation\label{natural}
\Theta(\gamma\cdot W)=\gamma \Theta(W)\gamma^{-1}
\EndEquation
for every $\gamma\in\Gamma$ and every saturated set $W\subset|K|$. 
\EndNumber
}

\Number\label{irving}
Now consider a discrete, torsion-free subgroup $\Gamma$ of ${\rm Isom}_+(\bf H^3)$ which is cocompact and hence purely loxodromic. Note that 
$\Gamma$ is an \iccg.
Note also that if $C$ is a maximal cyclic subgroup of $\Gamma$, all non-trivial elements of $C$ have the same translation axis, which will be called the {\it translation axis} of $C$.

Given a real number $\lambda >0$, we define ${\mathcal C}_{\lambda}(\Gamma) \subseteq {\mathcal C}(\Gamma)$ to be the set of maximal cyclic subgroups $ C$ of $\Gamma$ 
such that a (loxodromic) generator of $ C$ has translation length less than $\lambda$. 

For every non-trivial element $\gamma$ of  $\Gamma$, we denote by $Z_{\lambda}(\gamma)$  the set  $\{ \tP \in {\bf H}^3 :
d(\tP,\gamma \cdot \tP) < \lambda \}$.
For each  $C\in {\mathcal C}(\Gamma)$, 
we set $Z_{\lambda}(C)= \bigcup_{1 \not = \gamma \in {\mathcal C}} Z_{\lambda}(\gamma)$. 
For each ${ C} \in {\mathcal C}_{\lambda}(\Gamma)$, there is a
non-trivial element $\gamma \in { C}$ with $Z_{\lambda} ({ C})=
Z_{\lambda} (\gamma)$; furthermore,  $Z_{\lambda} ({ C})$ is a {\it hyperbolic cylinder} in the sense that it is a metric neighborhood of a line (the translation axis of $C$) with strictly positive radius. Note that any hyperbolic cylinder is a convex subset of $\HH^3$. 

For any $\lambda >0$ we set  $\mathcal{Z}_{\lambda}(\Gamma)= (Z_{\lambda}(C))_{C\in \calc_\lambda(\Gamma)}$. 
The discreteness of $\Gamma$ implies that the family $({Z}_{\lambda}
(\gamma))_{1\neq\gamma\in\Gamma}$ is locally finite, and hence the family ${\mathcal Z}_{\lambda}
(\Gamma)$ is locally finite as well. 
\EndNumber

\Number \label{GammaNerve}
Suppose that $\Gamma$ is a  cocompact,
 discrete, torsion-free subgroup of ${\rm Isom}_+(\bf H^3)$, and that $\lambda$ is a positive real number with the property that the family ${\mathcal Z}_{\lambda}(\Gamma)$ covers ${\bf H}^3$. 
Set $K=K_{{\mathcal Z}_{\lambda}(\Gamma)}$. Recall from \ref{first nerve} that the nerve $K$ comes equipped with a bijection $C\to v_C$ from $\calc_\lambda(\Gamma)$, the index set of the covering ${\mathcal Z}_{\lambda}(\Gamma)$ of $\HH^3$, to the vertex set $K^{(0)}$ of $K$.   If we denote the inverse bijection by $v\mapsto C_v$, then {$(C_v)_{v\in K^{(0)}}$ is a
$\Gamma$-labeling of $K$ (so that $(K,(C_v)_{v\in K^{(0)}})$ is a
$\Gamma$-labeled complex).} The labeling $(C_v)_{v\in K^{(0)}}$ will be referred to as {\it the canonical labeling} of $K=K_{{\mathcal Z}_{\lambda}(\Gamma)}$.

Since $C\to v_C$ is a bijection, we may define an action of $\Gamma$ on $K^{(0)}$ by setting $\gamma\cdot v_C=v_{\gamma C\gamma^{-1}}$ for each $\gamma\in\Gamma$ and each $C\in\calc_\lambda(\Gamma)$. Now suppose that $v_{C_0},\ldots,v_{C_d}$ are the vertices of a simplex of $K$. By the definition of the nerve $K$, we have $Z_\lambda(C_0)\cap\cdots\cap Z_\lambda(C_d)\ne\emptyset$. 
But for $i=0,\ldots,d$ and for any $\gamma\in\Gamma$, we
have $\gamma(Z_\lambda(C_i))=Z_\lambda(\gamma C_i\gamma^{-1})$. 
Hence $Z_\lambda(\gamma C_0\gamma^{-1})\cap\cdots\cap Z_\lambda(\gamma C_d\gamma^{-1})=
\gamma(Z_\lambda(C_0)\cap\cdots\cap Z_\lambda(C_d))\ne\emptyset$, which implies that $\gamma\cdot v_{C_0},\ldots,\gamma\cdot v_{C_d}$ 
are the vertices of a simplex of $K$. It follows that the action of $\Gamma$ on $K^{(0)}$ extends to a simplicial action on $K${, which  will be referred to as the canonical action of $\Gamma$ on 
$K=K_{{\mathcal Z}_{\lambda}(\Gamma)}$.}

Since we have $\gamma\cdot v_C=v_{\gamma C\gamma^{-1}}$ for each $\gamma\in\Gamma$ and each $C\in\calc_\lambda(\Gamma)$, and since $v\to C_v$ is the inverse bijection to $C\to v_C$, we have
$C_{\gamma\cdot v}
= \gamma C_v \gamma ^{-1}$ for each $\gamma\in\Gamma$ and each $v\in K^{(0)}$. Thus the canonical action of $\Gamma$ on $K$ is a labeling-compatible action. 
\EndNumber

\Proposition\label{fin dim}
Let $\Gamma$ be a cocompact, discrete, torsion-free subgroup of ${\rm Isom}_+(\bf H^3)$. Let $\lambda>0$ be given, and suppose that the family ${\mathcal Z}_{\lambda}(\Gamma)$ covers ${\bf H}^3$. Then  $K_{{\mathcal Z}_{\lambda}(\Gamma)}$  is finite-dimensional.
\EndProposition

\Proof
Let $D$ be a fundamental domain for the action of $\Gamma$ on $\HH^3$. Since $D$ is precompact by the cocompactness of $\Gamma$, and $\calz_\lambda(\Gamma)=(Z_{\lambda}(C))_{C\in C_\lambda(\Gamma)}$ is locally finite, the set $\calc^0:=\{C\in C_\lambda(\Gamma): Z_{\lambda}(C)\cap D\ne\emptyset\}$ is finite. {Let $N$ denote the cardinality of $\calc^0$. }

For each ${{\tP}}\in\HH^3$, set $\calc^{{\tP}}=\{C\in C_\lambda(\Gamma): {{\tP}}\in Z_{\lambda}(C)$. For any ${{\tP}}\in D$ we have $\calc^{{\tP}}\subset\calc^0$, so that {$\calc^{{\tP}}$ has cardinality at most $N$.
}
 For any point ${{\tP}}\in\HH^3$ there is an element $\delta$ of $\Gamma$ such that ${{\tP}}':=\delta\cdot {{\tP}}\in D$; the inner automorphism $\gamma\mapsto\delta\gamma\delta^{-1}$ of $\Gamma$ then carries $\calc^{{{\tP}}'}$ onto $\calc^{{\tP}}$, so that {$\calc^{{\tP}}$ has cardinality at most $ N$.
}
 Since {this is true
}
 for every ${{\tP}}\in\HH^3$, we have $\dim
K_{{\mathcal Z}_{\lambda}(\Gamma))}
\le N<\infty$.
\EndProof

\Number \label{Sp} {Suppose  that $\Gamma$ is a cocompact, discrete, torsion-free subgroup of ${\rm Isom}_+(\bf H^3)$. } Given {any} real number $\lambda>0$ and {any} point $\tP \in{\bf H}^3$, we will denote by
$\tS_\lambda(\tP,\Gamma)$
the set of all elements $ C$ of $\mathcal C(\Gamma)$ such that for some $\gamma  \in { C}-\{1\}$ we have $d(\tP,\gamma \cdot \tP) < \lambda $. Note that $\tS_\lambda(\tP,\Gamma)\subset\calc_\lambda(\Gamma)$, and that for any
$C\in\calc_\lambda(\Gamma)$ we have
  $C\in \tS_{\lambda}(\tP,\Gamma)$ if and only if $\tP\in Z_\lambda(C)$. 

Now suppose that $M$ is a closed orientable hyperbolic $3$-manifold (so that $\pi_1(M)$ is an \iccg). For any real number $\lambda>0$ and any point $P\in M$, we will denote by $S_\lambda(P)\subset\pi_1(M,P)$ the set of all elements of $C\in\calc(\pi_1(M,P))$ such that some non-trivial element of $C$ is represented by a loop of length less than $\lambda$. Note that {if}
$M=\HH^3/\Gamma$, where $\Gamma\le {\rm Isom}_+(\bf H^3)$ is discrete, torsion-free  and cocompact, and if $P$ is the image of a point $\tP\in\HH^3$ under the quotient map, then {$\Gamma$, regarded as the covering group of the simply-connected based covering $(\HH^3,\tP)$ of $(M,P)$, is canonically isomorphic to $\pi_1(M,P)$; and
 the canonical isomorphism 
carries $\tS_\lambda(\tP,\Gamma)\subset\Gamma$ onto $S_\lambda(P)\subset\pi_1(M,P)$.}
\EndNumber

\section{{The} Contractibility Lemma}\label{contractibility section}

For the reader's convenience we recall the following proposition from \cite{Guzman}.

\Proposition \label{homoequiv}
Suppose $L$ is a simplicial complex and $\sigma$ a simplex of $L$ such that ${\rm
  link}_L(\sigma)$ is contractible. If $X$ is a saturated subset of $|L|$  that contains all the simplices of $L$ for which $\sigma$ is a face, then the inclusion map $X - \sigma \hookrightarrow X$ is a homotopy equivalence. \EndProposition
\Proof This is \cite[Proposition 3.7]{Guzman}.
\EndProof

\Proposition \label{better homoequiv}
Let $L_0$ be  a subcomplex of a finite-dimensional simplicial complex $L$.  Suppose that for every simplex $\sigma$  of $L_0$, the complex ${\rm
  link}_L(\sigma)$ is contractible. Then the inclusion map
 $|L|-|L_0| \hookrightarrow |L|$ is a homotopy equivalence. \EndProposition
\Proof Let $n$ denote the dimension of $L$. For each $k$ with $0\le k\le n+1$, let $T_k$ denote the set of  all (open) simplices of $L$ having dimension at least $k$ (so that $T_{n+1}=\emptyset$ and $T_0=L$). Set  $W_k=(|L|-|L_0|)\cup\bigcup_{\sigma\in T_k}\sigma\subset|L|$. We have $|L|-|L_0|=
W_{n+1}\subset W_n\subset\cdots\subset W_0=|L|$. Thus it suffices to show that $W_{k+1}\hookrightarrow W_{k}$ is a homotopy equivalence for each $k$ with $0\le k\le n$.

Let us fix an index $k$ with $0\le k\le n$. The set $W_k-W_{k+1}$ is a union of (open) $k$-simplices. Let $\caly$ denote the set of all saturated subsets of $W_k$ having the form $W_{k+1}\cup F$, where $F$ is a finite union of $k$-simplices contained in $W_k-W_{k+1}$. Then $\caly$ may be regarded as a directed system in which the maps are inclusions, and $W_k$ is canonically identified with the direct limit of this system. Hence in order to show that $W_{k+1}\hookrightarrow W_{k}$ is a homotopy equivalence, we need only show that if $Y$ and $Y'$ are elements of $\caly$ and $Y$ is properly contained in $Y'$, then the inclusion map
$f:Y\hookrightarrow Y'$ is a homotopy equivalence. We may write $Y'-Y$ as a union of distinct $k$-simplices $\sigma_1,\ldots,\sigma_m$, where $m\ge1$. For each $i$ with $1\le i\le m$, set $X_i=Y\cup \sigma_1\cup\cdots\sigma_i$, and set $X_0=Y$. For $i=1,\ldots,m$, let $f_i:X_{i-1}\to X_i$ denote the inclusion map. Then $f=f_m\circ\cdots\circ f_1$.

For any index $i\in\{1,\ldots,m\}$, we have $X_{i-1}=X_i-\sigma_i$. Now $\sigma_i$ is disjoint from $Y$, and since $ Y\in\caly$ we have $Y\supset W_{k+1}\supset|L|-|L_0|$. Hence $\sigma_i\in L_0$. According to the hypothesis this implies that $\link_L(\sigma_i)$ is contractible. 
Note also that if $\tau$ is a simplex of $L$ having $\sigma_i$ as a face, then either $\tau=\sigma_i\subset X_i$, or $\dim\tau\ge1+\dim\sigma_i=k+1$. In the latter case we have $\tau\subset W_{k+1}$ (since $\tau\in T_{k+1}$), and $W_{k+1}\subset Y$ since $Y\in\caly$. Since $Y\subset X_i$ it follows that $\tau\subset X_i$ in this case. Thus every simplex of $L$ having $\sigma_i$ as a face is contained in $X_i$. We 
have now shown that the hypotheses of Proposition \ref{homoequiv} hold with $X=X_i$ and $\sigma=\sigma_i$. Hence the inclusion map $f_i$ from $X_{i-1}=X_i-\sigma_i$ to $X_i$ is a homotopy equivalence. Since this is true for $i=1,\ldots,m$, it follows that
$f=f_m\circ\cdots\circ f_1$ is a homotopy equivalence, as required.

\EndProof

The following lemma, which was mentioned in the introduction, is an analogue of \cite[Lemma 3.8]{Guzman}  in which ``minimum enveloping rank'' is replaced by internal rank; the argument is not affected in a radical way by this change.

\Lemma\label{ContractibilityLemma} Suppose that $\Gamma$ is a cocompact, torsion-free, discrete subgroup of ${\rm Isom}_+(\bf H^3)$. Let $\lambda$ be a positive real number, and let $m$ be a strictly positive integer.  Suppose that for every $\tP \in {\bf H}^3$, we have ${\rm IR}({\tS_\lambda(\tP,\Gamma)}) \ge m$. Set ${\mathcal C}= {\mathcal C}_{\lambda}(\Gamma)$. Then:
\begin{enumerate}
\item $\calz:={\mathcal Z}_{\lambda}(\Gamma)$
is a cover of ${\bf H}^3$ by cylinders. 
\item  The nerve $K:=K_\calz$ is contractible.
\item If $K$ is equipped with the canonical $\Gamma$-labeling described in \ref{GammaNerve}, so that {$\shmel{K}{m-1}$} is defined by \ref{more Theta}, then ${\rm link}_K(\sigma)$ is contractible for every simplex $\sigma$ of {$\shmel{K}{m-1}$.}
\item The space
{$|K|-|\shmel{K}{m-1}|$}
is contractible. 
\end{enumerate}
\EndLemma

\Proof 
By hypothesis we have ${\rm IR}({\tS_\lambda(\tP,\Gamma)}) \ge m \ge 1$ for every point ${\tP} \in {\bf H}^3$. Hence for any ${\tP} \in {\bf H}^3$, we have  {$\tS_\lambda(\tP,\Gamma)\ne\emptyset$, which by the discussion in \ref{Sp} implies that
$\tP \in Z_{\lambda}(C)$} for some $C \in {\mathcal C}$. Therefore, the family ${\mathcal Z}=(Z_{\lambda}(C))_{C \in {\mathcal C}}$ covers ${\bf H}^3$. This proves (1).

For every $C \in {\mathcal C}$, the cylinder $Z_{\lambda}(C)$ is open and convex. It follows that every nonempty, finite intersection of sets in the open covering ${\mathcal Z}$ is convex and therefore contractible. Thus, the Borsuk Nerve Theorem (see \ref{first nerve}) applies, and gives that $|K|$ is homotopy-equivalent to ${\bf H}^3$, and therefore contractible. This proves (2).

We now turn to the proof of (3).
Suppose that $\sigma$ is an open simplex of 
{$\shmel{K}{m-1}$}, which by definition {means} that the internal rank of $\sigma$ is at most $m-1$. Recall from \ref{first nerve} that the nerve $K$ comes equipped with a 
 bijection $C\mapsto v_C$ from the index set  ${\mathcal C}$ of $\calz$ to the vertex set $K^{(0)}$ of $K$, and from \ref{more Theta} that  
$I_\sigma= \{ C_v : v \in \sigmabar^{(0)}\} $. 
According to the definition of the canonical labeling of $K$ given in \ref{GammaNerve}, the assignment
$v\mapsto C_v$ is the inverse bijection to $C\mapsto v_C$.
Hence we have  
 $I_{\sigma} = \{ C \in {\mathcal C} : v_C \in \sigmabar^{(0)} \}$.

{Let $C_0, \dots, C_l$ denote the elements of $I_{\sigma}$, so that $v_{C_0 },\ldots,v_{C_l}$ are the vertices of $\sigma$.}
Set ${\mathcal U}_{\sigma} = Z_{\lambda}(C_0)  \cap \cdots \cap Z_{\lambda}(C_l)$.  Since $\sigma$ is in particular a simplex of the nerve $K$ of $\calz$, we have $\calu_\sigma\ne\emptyset$. Furthermore, the open set $\calu_\sigma\subset\HH^3$ is convex since it is an intersection of cylinders.


Set $J_{\sigma} =
\{ C\in {\mathcal C} - I_{\sigma} : Z_{\lambda}({C}) \cap {\mathcal U}_{\sigma} \neq \emptyset \}$. Set { $V_{{C},\sigma} = Z_{\lambda}({C}) \cap {\mathcal U}_{\sigma} $
 for each $ {C} \in J_{\sigma}$,}  and ${\mathcal V}_{\sigma}=(V_{{C},\sigma})_{{C} \in J_{\sigma}}$.


We claim:
\Claim \label{Claim1} ${\mathcal V}_{\sigma}$ is a cover for ${\mathcal U}_{\sigma}$. \EndClaim

To prove \ref{Claim1}, we suppose that ${\mathcal V}_{\sigma}$ is \emph{not} a cover for ${\mathcal U}_{\sigma}$. Then there exists a point ${\tP}$ of ${\mathcal U}_{\sigma}$
such that ${\tP} \not\in Z_{\lambda}(C)$ for any $C \in {\calc}-I_{\sigma}$. {But since ${\tP} \in {\mathcal U}_{\sigma}$, we have ${\tP} \in Z_{\lambda}(C)$ for any $C \in I_{\sigma}$. Thus for every $C\in\calc$, we have ${\tP} \in Z_{\lambda}(C)$ if and only if $C \in I_{\sigma}$. 
In view of the discussion in \ref{Sp}, this says that 
${\tS_\lambda(\tP,\Gamma)} = I_{\sigma}$.}

By definition we have ${\rm IR}(\sigma)={\rm IR}(I_{\sigma})$. We therefore have ${\rm IR}(\sigma) = {\rm IR}({\tS_\lambda(\tP,\Gamma)})\ge m$, where the last inequality {follows from the} hypothesis. But since $\sigma$ belongs to the subcomplex 
{$\shmel{K}{m-1}$},
 we have ${\rm IR} ({\tS_\lambda(\tP,\Gamma)})\le m-1$, a contradiction. Therefore, ${\mathcal V}_{\sigma}$ covers ${\mathcal U}_{\sigma}$, and \ref{Claim1} is proved.

According to \ref{first nerve}, the nerve {$K_{\calv_\sigma}$} comes equipped with a bijection from the index set $J_\sigma$ of the covering ${\mathcal V}_{\sigma}$ to the vertex set ${K_{\calv_\sigma}}^{(0)}$ of ${K_{\calv_\sigma}}$.  We will denote this bijection by $C\mapsto w_C$.


Since $J_\sigma\subset\calc$, and since $C\mapsto v_C$ is a bijection from the index set  ${\mathcal C}$ of $\calz_\lambda$ to the vertex set $K^{(0)}$ of $K$,
 there is a well-defined injection ${f_\sigma}^{(0)}:{K_{\calv_\sigma}}^{(0)}\to K^{(0)}$ given by ${f_\sigma}^{(0)}(w_C)=v_C$.

{
(Note that the indexing of the family $\calv_\sigma$ by $J_\sigma$ is needed to guarantee that the map $f_\sigma^{(0)}$ is well-defined. On the other hand, 
the indexed family $\calv_\sigma$ may fail to be ``faithfully indexed'' in the sense explained in \ref{first nerve}.
Specifically, if we have $Z_\lambda(C)\supset \calu_\sigma$ and $Z_\lambda(C')\supset \calu_\sigma$ for some $C,C'\in\calc$, then $V_{C,\sigma}=V_{C',\sigma}=\calu_\sigma$.} 
It is for this reason that we need to consider nerves of non-faithfully indexed covers, as mentioned in \ref{first nerve}.)

We claim:

\Claim\label{numbah one}The map ${f_\sigma}^{(0)}$ extends to an injective  simplicial map ${f_\sigma}:{K_{\calv_\sigma}}\to K$.
\EndClaim

To prove \ref{numbah one}, first note that if ${f_\sigma}^{(0)}$ admits a simplicial extension ${f_\sigma}:{K_{\calv_\sigma}}\to K$, then the injectivity of ${f_\sigma}^{(0)}$ implies that ${f_\sigma}$ is injective.  To prove the existence of the simplicial extension, we must show that if $w_0,\ldots,w_d$ are vertices of a simplex of ${K_{\calv_\sigma}}$, then ${f_\sigma}(w_0),\ldots,{f_\sigma}(w_d)$ are vertices of a simplex of $K$. For $i=0,\ldots,d$ let us write $w_i=w_{C'_i}$, where $C'_i\in J_\sigma$. We then have ${f_\sigma}(w_i)=v_{C'_i}$ for $i=0,\ldots,d$. Since ${K_{\calv_\sigma}}$ is the nerve of ${\mathcal V}_{\sigma}=(V_{j,\sigma})_{j \in J_{\sigma}}$, the condition that $w_{C'_0},\ldots,w_{C'_d}$ are vertices of a simplex of ${K_{\calv_\sigma}}$ means that $V_{C'_0,\sigma}\cap\cdots\cap V_{C'_d,\sigma}\ne\emptyset$. But for $i=0,\ldots,d$ we have $V_{C'_i,\sigma}=Z_{\lambda}(C'_i) \cap {\mathcal U}_{\sigma}\subset Z_{\lambda}(C'_i)$, so that in particular $Z_\lambda(C'_0)\cap\cdots\cap Z_\lambda(C'_d)\ne\emptyset$. This means that $v_{C'_0},\ldots,v_{C'_d}$ are the vertices of a simplex of the nerve $K$ of 
$\mathcal{Z}= (Z_{\lambda}(C))_{C\in C'_\lambda(\Gamma)}$, and \ref{numbah one} is established.

Defining ${f_\sigma}$ as in \ref{numbah one}, we now claim:

\Claim\label{it's the link}
We have ${f_\sigma}({K_{\calv_\sigma}})={\rm link}_K(\sigma)$.
\EndClaim

To prove \ref{it's the link}, first consider an arbitrary simplex $\tau$ of ${K_{\calv_\sigma}}$. Let $w_0,\ldots,w_d$ denote the vertices of $\tau$. For $i=0,\ldots,d$ let us write $w_i=w_{C'_i}$, where $C'_i\in J_\sigma$. We then have ${f_\sigma}(w_i)=v_{C'_i}$ for $i=0,\ldots,d$. {By the definition of the nerve ${K_{\calv_\sigma}}$
 of ${\mathcal V}_{\sigma}=(V_{j,\sigma})_{j \in J_{\sigma}}$,}
the condition that $w_{C'_0},\ldots,w_{C'_d}$ are vertices of a simplex of ${K_{\calv_\sigma}}$ means that $V_{C'_0,\sigma}\cap\cdots\cap V_{C'_d,\sigma}\ne\emptyset$. But for $i=0,\ldots,d$ we have $V_{C'_i,\sigma}=Z_{\lambda}(C'_i) \cap {\mathcal U}_{\sigma}$, so that $Z_\lambda(C'_0)\cap\cdots\cap Z_\lambda(C'_d)\cap\calu_\sigma\ne\emptyset$. Since ${\mathcal U}_{\sigma} = Z_{\lambda}(C_0) 
\cap \cdots \cap Z_{\lambda}(C_l)$, this means that $Z_\lambda(C'_0)\cap\cdots\cap Z_\lambda(C'_d)\cap
Z_{\lambda}(C_0) 
 \cap \cdots \cap Z_{\lambda}(C_l)
\ne\emptyset$. Hence ${f_\sigma}(w_0),\ldots,{f_\sigma}(w_d),v_0,\ldots,v_l$ are the vertices of a simplex $\varphi$ of $K$; that is, ${f_\sigma}(\tau)$ and $\sigma$ are faces of $\varphi$. Therefore ${f_\sigma}(\tau)$ is a simplex of the closed star 
$\mathop{\rm Cl}\mathop{\rm St}_K(\sigma)$.  On the other hand, since $\tau$ is a simplex of ${K_{\calv_\sigma}}$, we have $C'_i\in J_\sigma\subset \calc-I_\sigma$ for $i=0,\ldots,d$; hence ${f_\sigma}(w_i)$ is not a vertex of $\sigma$ for $i=0,\ldots,d$. This means that the  simplex ${f_\sigma}(\tau)\in\mathop{\rm Cl}\mathop{\rm St}_K(\sigma)$  has no vertices in common with  $\sigma$, and hence ${f_\sigma}(\tau)$ is a simplex of ${\rm link}_K(\sigma)$. Thus we have shown that ${f_\sigma}({K_{\calv_\sigma}})\subset{\rm link}_K(\sigma)$.

To establish the reverse inclusion, consider an arbitrary  simplex $\psi$ of ${\rm link}_K(\sigma)$. Let $z_0,\ldots,z_d$ denote the vertices of $\psi$.  For $i=0,\ldots,d$ let us write $z_i=v_{C'_i}$, where $C'_i\in I_\sigma$. Since
$\psi\in{\rm link}_K(\sigma)$, the simplices $\sigma$ and $\psi$ are faces of a simplex $\varphi$ of $K$. This means that $z_0,\ldots,z_d,v_0,\ldots,v_l$ are  vertices of  $\varphi$; thus we have $(Z_\lambda(C'_0)\cap\cdots\cap Z_\lambda(C'_d))\cap
(Z_{\lambda}(C_0) \cap \cdots \cap Z_{\lambda}(C_l))
\ne\emptyset$, i.e.
\Equation\label{i.e.}
Z_\lambda(C'_0)\cap\cdots\cap Z_\lambda(C'_d)\cap
\calu_\sigma
\ne\emptyset.
\EndEquation

It follows from (\ref{i.e.}) that 
$Z_\lambda(C'_i)\cap\calu_\sigma\ne\emptyset$ for $i=0,\ldots,d$. On the other hand, since 
$\psi\in{\rm link}_K(\sigma)$, the simplices $\sigma$ and $\psi$ have no vertices in common, and hence $C'_i\in\calc-I_\sigma$ for $i=0,\ldots,d$. In view of the definition of $J_\sigma$, it now follows that $C'_i\in J_\sigma$ for $i=0,\ldots,d$. Thus $w_{C_0'},\ldots,w_{C_d'}$ are vertices of ${K_{\calv_\sigma}}$, and ${f_\sigma}^{(0)}(w_{C_i'})=z_i$ for $i=0,\ldots,d$. But by definition we have  $V_{C_i',\sigma} = Z_{\lambda}(C_i') \cap {\mathcal U}_{\sigma} $
 for $i=0,\ldots,d$, so that (\ref{i.e.}) gives  
$V_{C_0',\sigma}\cap\cdots\cap V_{C_d',\sigma} \ne\emptyset$. Hence 
$w_{C_0'},\ldots,w_{C_d'}$ are the vertices of a simplex $\tau$ of ${K_{\calv_\sigma}}$. Since  ${f_\sigma}^{(0)}(w_{C_i'})=z_i$ for each $i$, we have ${f_\sigma}(\tau)=\psi$. This shows that ${\rm link}_K(\sigma) \subset {f_\sigma}({K_{\calv_\sigma}})$, and completes the proof of (\ref{i.e.}).

We are now in a position to show that ${\rm link}_K(\sigma)$ is contractible, which will establish Assertion (3).
It follows from \ref{numbah one} and \ref{it's the link} that ${f_\sigma}$ may be regarded as a simplicial isomorphism of  ${K_{\calv_\sigma}}$ onto ${\rm link}_K(\sigma)$. In particular, ${\rm link}_K(\sigma)$ is homeomorphic to ${K_{\calv_\sigma}}$. But  ${K_{\calv_\sigma}}$ is the nerve of the covering $\calv_\sigma=(V_{C,\sigma})_{C \in J_{\sigma}}$
 of $\calu_\sigma$. We have observed that $\calu_\sigma$ is a nonempty convex subset of $\HH^3$; hence it is contractible. For each $C \in J_{\sigma}$, the set $V_{C,\sigma}=Z_{\lambda}(C) \cap {\mathcal U}_{\sigma} $ is the intersection of the convex set $\calu_\sigma$ with a hyperbolic cylinder, and is therefore convex. Hence every nonempty finite intersection of sets in the family $\calv_\sigma$ is convex and therefore contractible. It now follows from the Borsuk Nerve Theorem that ${K_{\calv_\sigma}}$ is contractible. Hence ${\rm link}_K(\sigma)$ is contractible, and (3) is established.




To prove (4), note that according to Lemma \ref{fin dim} and Assertion (1), the complex $K$ is finite-dimensional. Since by Assertion (3) every simplex of $K_{({m-1})}'$ has a 
contractible link in $K$, we may apply Proposition \ref{better homoequiv}, with $K$ and $K_{({m-1})}'$ playing the respective roles of $L$ and $L_0$, to deduce that the inclusion map $K_{({m-1})}'\hookrightarrow K$ is a homotopy equivalence. Since $K$ is contractible by Assertion (2), it now follows that $K_{({m-1})}'$ is contractible, and (4) is established.
\EndProof


\section{{The} Local Rank Lemma}\label{rank section}
\Definition \label{local rank} A group $G$ is said to have {\it local rank} $\le k$,  where $k$ is a non-negative integer, if every finitely generated subgroup of $G$ is
contained in a subgroup of $G$ which has rank less than or equal
to $k$. The local rank of $G$ is the smallest integer $k\ge0$ with this
property. If there does not exist such a $k$ then we define the local
rank of $G$ to be $\infty$. Notice that if $G$ is
finitely generated, its local rank is simply its rank. \EndDefinition

The following lemma, which was discussed in the introduction, is an analogue of  \cite[Proposition 4.4]{cs} and \cite[Proposition 4.3]{Guzman}, in which ``minimum enveloping rank'' is replaced by internal rank; as we pointed out in the introduction,  this change has a major effect on the argument, in that it obviates the need for results comparing ranks of joins and intersections in free groups.

%

\Lemma \label{RankLemma} 
{Let $k$ be a positive integer.} 
Suppose that $\Gamma$ is a $k$-free group, and that
$({L}, (C_v)_{{v\in L^{(0)}}})$ a $\Gamma$-labeled complex. Fix an integer $r$ with {$0< r\le k-1$}. Let $W$ be a subset of $|{L}|$ which is saturated and connected, and has the property that {every simplex}  in $W$ {has} internal rank exactly $r$. Then the local rank of $\Theta(W)$ is at most $ r$ (where $\Theta(W)$ is defined by \ref{more Theta}).
\EndLemma




\Proof[Proof of Lemma \ref{RankLemma}] In view of the definition of local rank, we must show that every finitely generated subgroup of
$\Theta(W)$ is contained in a finitely generated subgroup of $\Theta(W)$ of rank $\le r$. So suppose that $E \le \Theta(W)$ is a finitely
generated subgroup of $\Theta(W)$. Since $W$ is connected, we have $E\le \Theta(V)$ for some
saturated subset $V$ of $W$ that is connected and contains only finitely many open
simplices of $W$. 
In particular, $\Theta(V)$ is finitely generated.

Since $V$ is connected, we
may {index} the  open simplices of $V$ as $\sigma_0,\dots,\sigma_m$, {where $m$ is a non-negative integer,} {in such a way that for every index} $i$ with $0 < i \le m$, there is an
index $l$ with $0 \le l < i$ such that {either} $\sigma_l$ is a proper face
of $\sigma_i$ or $\sigma_i$ is a proper face of $\sigma_l$. Define
the saturated subset $V_i = \sigma_0 \cup \dots \cup \sigma_i$ for
$0\le i\le m$; by induction on $i$, we will show $\rank
\Theta(V_i) \le  r$. { For the base case, note that by Definition \ref{more Theta}, and the hypothesis of the present lemma, we have 
$\rank\Theta(V_0)= \rank\Theta(\sigma_0)\le
\max \{\rank\Theta(\tau) : \tau \le \sigma_0 \}=
{\rm IR} (\sigma_0) = r$.}

For the induction
step we assume $\rank\Theta(V_{i-1}) \le r$. We want to show that
$\rank\Theta(V_i) \le r$. By the way we have ordered the list of simplices in $V$, there is an index
$l$ with $0 \le l < i$ such that $\sigma_l$ is a proper face of
$\sigma_i$ or $\sigma_i$ is a proper face of $\sigma_l$.

First consider the case in which $\sigma_i$ is a proper face
of $\sigma_l$. Then the set of vertices of simplices in $V_i$ is the same as 
the set of vertices of simplices in $V_{i-1}${; hence in the notation of Definition \ref{more Theta}, we have $I_{V_{i}}=I_{V_{i-1}}$, so} that
$\Theta(V_i) {=\langle I_{V_i}\rangle=\langle I_{V_{i-1}}\rangle=  \Theta(V_{i-1})}$.
By the induction hypothesis we have  ${\rm rank}
\;\Theta(V_{i-1}) \le r$, and so ${\rm rank} \;\Theta(V_i) \le r$ as
required.

 Next, consider the case in which $\sigma_l$ is a proper
face of $\sigma_i$. We have ${\rm rank} \;\Theta(V_{i-1}) \le r$ by the induction hypothesis and ${\rm IR}
(\sigma_i) = {\rm IR} (\sigma_l) = r$ by the hypothesis of the lemma. We need to show $\Theta
(V_i) = \Theta(V_{i-1}) \vee \Theta(\sigma_i)$ has rank $\le r$.

Let $ v_{1},...,v_{n} $ be the vertices of $\sigma_i$ which
are not vertices of $\sigma_l$. Define $X_j= C_{v_{j}}$ for $j=1,\ldots, n$.
Set $X_0 =\Theta(V_{i-1})$, and let $Y_m =\langle X_0, X_1\dots X_m \rangle$ for $m=0,\ldots, n$. Then $Y_n=\Theta(V_{i})$. We will assume that $\rank \Theta(V_i)=\rank Y_n>r$ and obtain a contradiction. The assumption implies that there is a smallest index $t\in\{0,\ldots,n\}$ such that $\rank Y_t>r$. Since $\rank Y_0=\rank X_0=\rank\Theta(V_{i-1})\le r$, we have $t\ge1$. The minimality of $t$ implies that $\rank Y_{t-1}\le r$. We have $Y_t=\langle Y_{t-1},X_t\rangle$, so that $\rank Y_t\le 1+\rank Y_{t-1}$. Since $\rank Y_t>r\ge\rank Y_{t-1}$, we must have $\rank Y_{t-1}=r$ and $\rank Y_{t}=r+1$. On the other hand, since $r+1\le k$, and $\Gamma$ is $k$-free, the subgroup $Y_t$ is free. Fix a free basis $u_1,\ldots,u_r$ of $Y_{t-1}\le Y_t$, and fix a generator $x$ of $X_t$. Then $\{u_1,\ldots,u_r,x\}$ is a generating set of cardinality $r+1$ for the rank-$(r+1)$ free group $Y_t$, and is therefore a free basis. This shows that
$Y_t$  a free product: $Y_t = Y_{t-1} \ast X_t$.


Because ${\rm IR}(\sigma_l) = r$, by definition there exists a face  $\tau$  of $\sigma_l$ with ${\rm rank}\;\Theta(\tau) = r$. Let $\sigma '$ denote the face of $\sigma_i$ which is spanned by the face $\tau$ and the vertex $v_t$; thus 
$\Theta(\sigma ') = \langle \Theta(\tau) , C_{v_t}\rangle$. But since $\tau\le\sigma_l$, we have  $\Theta(\tau)\le\Theta(\sigma_l)\le X_0\le Y_{t-1}$. Since $Y_t = Y_{t-1} \ast X_t$, it now follows that  $\Theta(\sigma ') = \Theta(\tau) \ast C_{v_t}$. Hence $\rank\Theta(\sigma')=\rank\Theta(\tau)+1=r+1$.
However, $\sigma '$ is a face of $\sigma_i$, and since $\sigma_i\subset W$ we have ${\rm IR}(\sigma_i)=r$ by hypothesis. In view of the definition of internal rank, it follows that $\rank\Theta(\sigma ') \le r$,  a contradiction. 
\EndProof

\section{Proof of the Main Theorem} \label{main proof section}

{As a preliminary to the proof of our main theorem, we review the following result from \cite{Guzman}:

\Lemma\label{guzman two-five}
Suppose that $k \ge 2$ is an integer, and that 
$\Gamma\le 
\Isom_+(\HH^3)$ is discrete, loxodromic, and $k$-free (and therefore torsion-free). If $C_1,\ldots,C_n$ are elements of $\calc_{\log(2k-1)}(\Gamma)$ such that  $Z_{\log (2k-1)}(C_1)\cap\cdots\cap Z_{\log (2k-1)}(C_n)\ne\emptyset$, then the rank of $\langle C_1,\ldots , C_n\rangle$ is at most $ k- 1$.
\EndLemma

\Proof
This is Lemma 2.5 of \cite{Guzman}.
\EndProof
}

\label{Final}
Let us restate Main Theorem \ref{Main Theorem} from the introduction in terms of the notion of internal rank, and the notation $S_\lambda (P)$, which were introduced in Section \ref{prelim section}:

\Theorem \label{MainTheoremRestated} Suppose {that $k\ge 3$ is an integer, and that} $M $ is a $k$-free, closed, orientable, hyperbolic 3-manifold.  Then there exists a point $P$ of $M$ such that \linebreak ${\rm IR}(S_{\log(2k-1)}(P)) \le k-3$.
\EndTheorem

\Proof
Let us write $M = {\bf H}^3/\Gamma$, where $\Gamma \le {\rm Isom}_+
({\bf H}^3)$ is discrete, cocompact, and $k$-free (and in particular torsion-free). According to the  discussion in \ref{Sp}, it suffices to show that there is a point ${\tP}\in\HH^3$ such that  ${\rm IR}(\tS_{\log(2k-1)}(\tP,\Gamma)) \le k-3$.

Arguing by contradiction, assume that for every point $\tP\in\HH^3$ we have ${\rm IR}(\tS_{\log(2k-1)}(\tP,\Gamma) \ge k-2$. According to Lemma \ref{ContractibilityLemma}, applied with $m=k-2$ {(which is strictly positive since $k\ge3$)}, it follows both that ${\mathcal Z}:=\calz_{\log{(2k-1)}} (\Gamma)$ is a cover of 
${\bf H}^3$, and that if we set $K=K_\calz$, then
{$|K|-|
\shmel{K}{k-3}|$}
 is contractible.

We claim:

\Claim\label{less than k}
For every simplex $\sigma \subset
{|K|-|
\shmel{K}{k-3}|}$,
we have
${\rm IR}(\Theta(\sigma)) \le  k-1$.
\EndClaim

To prove \ref{less than k}, first recall from \ref{more Theta} that  ${\rm IR}(\sigma)=\max_{\tau\le\sigma}\rank \Theta (\tau)$. Hence we need only show that if $\tau$ is any face of $\sigma$ then $\rank\Theta(\tau)\le k-1$. By definition, $\Theta(\tau)=\langle C_{v_1},\ldots,C_{v_n}\rangle$, where
$v_1,\ldots,v_n$ are the vertices of $\tau$. Since $\tau$ is a simplex of the nerve $K$ of $\calz$, we have $Z_{\log(2k-1)}(C_{v_1})\cap\ldots\cap Z_{\log(2k-1)}(C_{v_n})\ne\emptyset$. Since $\Gamma$ is $k$-free, it then follows from Lemma \ref{guzman two-five} that $\rank\langle C_{v_1},\ldots,C_{v_n}\rangle\le k-1$. Thus \ref{less than k} is established.

{It follows from the definition of 
{$
\shmel{K}{k-3}$}
(see \ref{more Theta})  that a simplex is contained in {$|K|-|
\shmel{K}{k-3}|$}
if and only if its internal rank is at least $k-2$. Combining this with \ref{less than k}, we deduce that the saturated subset 
{$|K|-|
\shmel{K}{k-3}|$}
is the union of simplices of $K$ of internal rank equal to $k-2$ or $k-1$. For $r=k-1,k-2$, we let $X_r$ denote the union of simplices of internal rank equal to $r$, so that 
{$|K|-|
\shmel{K}{k-3}|$}
is the set-theoretical disjoint union of $X_{k-2}$ and $X_{k-1}$.}

For each $r \in \{ k-2, k-1 \}$, let $\mathcal W_r$ denote the set of connected components of $X_r$. 
We will say that elements $W_{k-2}$ of $\calw_{k-2}$ and $W_{k-1}$ of $\calw_{k-1}$ are {\it adjacent} if  there  are simplices $\sigma_{k-2}\subset{ W}_{k-2}$ and
 $\sigma_{k-1}\subset{ W}_{k-1}$ such that either $\sigma_{k-2}<\sigma_{k-1}$ or $\sigma_{k-1}<\sigma_{k-2}$.

We construct an abstract bipartite graph ${\mathcal G}$ as follows: The vertex set $Z$ of ${\mathcal G}$ is a disjoint union $Z_{k-2}\,\Dot\cup\, Z_{k-1}$, where for $r=k-2,k-1$ the set $Z_r$ is a bijective copy of  $\mathcal W_{r}$ and is equipped with a bijection $W\mapsto s_W$ from $\calw_r$ to $Z_r$. An edge of $\calg$ is defined to be an unordered pair of the form $\{s_{W_{k-2}},s_{W_{k-1}}\}$, where $W_r$ is an element of $\calw_r$ for $r=k-2,k-1$, and $W_{k-2}$ and $W_{k-1}$ are adjacent. 

Let $T$ denote the geometric realization of $\calg$ (regarded as a simplicial complex of dimension at most $1$).
Because 
{$|K|-|
\shmel{K}{k-3}|$}
 is the disjoint union of the saturated subsets  $X_{k-2}$ and $ X_{k-1}$ of $|K|$, it follows from \cite[Lemma 5.12]{cs} that $|{T}|$ is a homotopy-retract of 
{$|K|-|
\shmel{K}{k-3}|$}.
Since 
{$|K|-|
\shmel{K}{k-3}|$}
is contractible, it follows that
${T} $ is a tree. 

According to \ref{GammaNerve},  the canonical action of
$\Gamma$ on $K$ is labeling-compatible. Hence for any $\gamma \in \Gamma$ and any simplex
$\tau$ of $K$, we may apply (\ref{natural}), with $\tau$ playing the role of $X$, to deduce that 
$\Theta(\gamma\cdot\tau)=\gamma\Theta(\tau)\gamma^{-1}$. In particular 
$\Theta(\gamma\cdot\tau)$ and $\Theta(\tau)$
have the same rank. Hence if $\sigma$ is any simplex of $K$, we have ${\rm IR}(\gamma\cdot\sigma)=\max_{\tau\le\sigma}\rank\Theta(\gamma\cdot\tau)=
\max_{\tau\le\sigma}\rank\Theta(\tau)={\rm IR}(\sigma)$.
Consequently, $X_{k-2}$ and $X_{k-1}$ are invariant under the action of $\Gamma$. 

Since the canonical action of $\Gamma$ on $K$ is simplicial, it defines a continuous action of $\Gamma$ on $|K|$; hence for $r=k-2,k-1$, the restricted action of $\Gamma$ on $X_r$ gives rise to an action on $\calw_r$. The simplicial nature of the action on $K$ also implies that if $W_{k-1}\in\calw_{k-1}$ and $W_{k-2}\in\calw_{k-2}$ are adjacent, then $\gamma\cdot W_{k-1}$ and $\gamma\cdot W_{k-2}$ are adjacent for every $\gamma\in\Gamma$. Hence the action of $\Gamma$ on $Z$ defined by $\gamma\cdot z_W=z_{\gamma\cdot W}$ extends to a simplicial action of $\Gamma$ on $\calg$, whose geometric realization is an action on $T$.

Note that since $X_{k-2}$ and $X_{k-1}$ are $\Gamma$-invariant, the sets $Z_{k-2}, Z_{k-1}\subset\calg$ are also $\Gamma$-invariant. Since each edge of $\calg$ has one vertex in $Z_{k-2}$ and one in $Z_{k-1}$, the action of $\Gamma$ on $T$ has no inversions.

We claim:
\Claim\label{have a stab}
The stabilizer in $\Gamma$ of each vertex of $T$ is locally free.
\EndClaim

To prove \ref{have a stab}, consider an arbitrary vertex $s$ of $\calg$, and write $s=s_W$ for some $W\in\calw_r$, where $r\in\{k-2,k-1\}$. If $\gamma$ is any element of the stabilizer $\Gamma_s\le\Gamma$, it follows from the definition of the action of $\Gamma$ on $\calg$ that $\gamma\cdot W=W$. On the other hand, applying (\ref{natural}), with  $W$ defined as above, we find that 
$\Theta(\gamma\cdot W)=\gamma\Theta(W)\gamma^{-1}$. Hence
$\Theta( W)=\gamma\Theta(W)\gamma^{-1}$. This shows that $\Gamma_s$ is contained in the normalizer of $\Theta(W)$, which we shall denote by $H$.


Since $W\in\calw_r$, the set $W$ is a (nonempty) connected, saturated subset of $|K|$, and every simplex contained in $W$ has internal rank $r$. We have $r<k$. By Lemma \ref{RankLemma}, applied with $L=K$ and with $W$ defined as above,  the local rank {$l$} of $\Theta(W)$  must be at most $r$. In particular ${l}<k$. 
Since $W$ is nonempty it must contain at least one simplex $\tau$, which must have at least one vertex $v_0$; we have $\Theta(W)\ge\Theta(\tau)\ge C_{v_0}$, and hence $\Theta(W)\ne\{1\}$. It follows that ${l}\ge1$.

According to \cite[Proposition 4.5]{cs}, if a $k$-free group $G$ has a normal subgroup $N$ of  local rank  ${l}$, where ${l}$ is an integer with $1\le {l}<k$, then the group $G$ is itself of local rank at most $ {l}$. 
We shall apply this result, taking $G=H$ and taking  $N=\Theta(W)$. Since $\Gamma$ is $k$-free, its subgroup $H$ is also $k$-free. The definition of $H$ implies that $\Theta(W)$ is a normal subgroup of $H$. We have seen that the local rank ${l}$ of $\Theta(W)$ satisfies $1\le {l}<k$. Thus all the hypotheses of  
\cite[Proposition 4.5]{cs} hold in this situation, and 
we deduce that the local rank of $H$ is at most ${l}\le k-1$. In particular, $H$ is locally free, being a subgroup of local rank $\le k$ of the $k$-free group $\Gamma$. Therefore $\Gamma_s\le H$ is also locally free. This completes the proof of \ref{have a stab}.

We have shown that $\Gamma\cong\pi_1(M)$ acts simplicially, without inversions, on the tree $T$, and by \ref{have a stab}, the stabilizer of every vertex in $\Gamma$ is locally free. But according to \cite[Lemma 5.13]{cs}, the fundamental group of a closed, orientable, aspherical $3$-manifold $M$ cannot admit a simplicial action, without inversions, on a tree in such a way that each vertex stabilizer is locally free. This  contradiction completes the proof.
\EndProof

\smallskip
\bibliographystyle{plain}

\end{document}